\documentclass[a4paper,12pt]{article}
\usepackage{amsmath}
\usepackage{amsfonts}
\usepackage{amsthm}
\usepackage{graphicx}

\newtheorem{lemma}{Lemma}[section]
\newtheorem{propos}[lemma]{Proposition}

\newtheorem{theorem}[lemma]{Theorem}
\newtheorem{cor}[lemma]{Corollary}

\newcommand{\rat}{\ensuremath{\mathbb{Q}}}

\newcommand{\axref}{\ensuremath{\mathit{Ref}}}
\newcommand{\ril}{\ensuremath{\mathit{Ril}}}
\newcommand{\pil}{\ensuremath{\mathit{Pil}}}

\newcommand{\bfi}{\ensuremath{\mathit{\bf i}}}
\newcommand{\bfj}{\ensuremath{\mathit{\bf j}}}
\newcommand{\bfk}{\ensuremath{\mathit{\bf k}}}

\begin{document}


\begin{center} {\Large {\bf Skew Meadows}}\\
\today 
   \\ \baselineskip 13pt{\ } \vskip 0.2in
J A Bergstra\footnote{Email: j.a.bergstra@uva.nl}
\\{\ }\\ Informatics Institute, University of Amsterdam,  \\ Science Park 403, 
		1098 SJ  Amsterdam,  The Netherlands
\\{\ }\\ Y Hirshfeld\footnote{Email: joram@post.tau.ac.il}
\\{\ }\\ Department of Mathematics, Tel Aviv University, \\Tel Aviv 69978, Israel
\\{\ }\\ J V Tucker\footnote{Email: j.v.tucker@swansea.ac.uk}
\\{\ }\\ Department of Computer Science, Swansea University, \\ Singleton Park,  Swansea, SA2 8PP, United Kingdom
\end{center}

\bigskip

\begin{abstract}
\noindent
A skew meadow is a non-commutative ring with an inverse operator satisfying two special equations and in 
which $0^{-1} = 0$.  All skew fields and products of skew fields can be viewed as skew meadows. Conversely, we give an embedding of non-trivial skew meadows into products of skew fields, from which a completeness result for the equational logic of skew fields is derived. The relationship between regularity conditions on rings and skew meadows is investigated.
\end{abstract}
\bigskip

\section{Introduction}
Meadows have been introduced in \cite{BT07} while the equations for meadows used below were improved in \cite{BHT07,BHT08}. Meadows may be viewed as a generalization of so-called zero totalized fields, being fields in which division is made total by setting $0^{-1}=0$. Thanks to a characterization theorem in \cite{BHT07}, however,  a {\em meadow} can also be defined as a commutative ring with unit equipped with a total unary operation $x^{-1}$, 
named inverse, that satisfies these two additional equations:
\begin{eqnarray}
(x^{-1})^{-1}& = & x \\
x \cdot (x \cdot x^{-1}) & = & x
\end{eqnarray}
The first equation is called \axref, the second equation is called \ril. 

Perhaps the clearest way to specify the class of meadows is as the smallest variety containing all zero-totalized fields. This is a matter of taste to some extent.

Commutative meadows provide an analysis of division which is more general than that of the classical theory of fields. Commutative meadows are total algebras in which $0^{-1} = 0$. We have used algebras with such zero totalized division in developing elementary algebraic specifications for several algebras of rational numbers in our previous paper \cite{BT07} and its companions
  \cite{B06, BT06}.
  
Several generalizations of meadows can be conceived. In \cite{BP08} signed meadows were defined as a generalization of ordered fields (with totalized division), and in \cite{BP08b} differential meadows are considered which may be viewed as a generalization of differential fields. 

The generalization of commutative meadows, as defined in \cite{BT07} and subsequently analyzed in \cite{BHT08}, to the non-commutative case is the subject of this paper.  As is always the case, the transition from commutative to non-commuutative rings is a delicate operation, leading to a ramification of properties. However, we are able isolate a number of concepts and prove nice generalizations of basic results, including the following Representation Theorem \ref{EMB}: 
\\
\\
{\bf Theorem}  \emph{An algebra is a non-trivial skew meadow if, and only if, it is isomorphic to a subalgebra of a product of zero totalized skew fields.}
\\
\\
Ring theory has several concepts, like Von Neumann regularity, that distinguish elements with properties similar to multiplicative inverses, but does not seem to have investigated the possible corresponding inverse operators. We establish the relationship between skew meadows and forms of von Neumann regularity conditions on non-commutative rings (Theorem \ref{SRRtoSkMd}). The equational nature of meadows is confirmed for skew meadows by this result (Theorem \ref{comp}): 
\\
\\
{\bf Theorem} \emph{An equation is valid in all skew fields with zero totalized division if, and only if, it is true in the variety of skew meadows and for that reason derivable from their equational axiomatisation.}
\\
\\
Investigations of fields with zero totalized division have some history in logic and computing: it is mentioned in \cite{Ho93} as a reasonable method to extend division to a total function, and it has been used in a more technical way in \cite{H1998}. This work may be viewed as belonging to universal algebra and equational logic, with some orientation towards computer science, in particular, to the theory of abstract data types. We refer to \cite{MT92} and \cite{W92} for computer science oriented introductions to universal algebra. Of course the paper is about noncommutative rings. We refer to \cite{McC64} and \cite{La01} for introductions to noncommutative rings. To some extent the results are intimately related to the theory of regular and inverse semigroups, because many of the arguments can be given without a reference to addition and subtraction \cite{Law98}. We will focus, however, on connections with the theory of noncommutative rings.

\section{Axioms for rings}\label{ANA}
We start with a listing of the axioms of a unital ring. The starting point is a signature $\Sigma_{RU}$ for rings with unit:
\newline

{\bf signature}  $\Sigma_{RU}$

{\bf sorts}     $ring$

{\bf operations}  

$0 \colon  \to ring$;

$1 \colon  \to ring$;

$ + \colon  ring \times ring \to ring$;

$ - \colon  ring \to ring$;

$ \cdot \colon  ring \times ring\to ring$

{\bf end}\\

The first set of axioms is that of a ring with $1$, which establishes the standard properties of $+$, $-$, and $\cdot$. 
\newline

{\bf equations} {\it $RU$}
\begin{eqnarray}
(x + y) + z & = & x + (y + z)    \\
x + y & = & y + x  \\
x + 0 & = & x  \\
x + (-x) & = & 0  \\
(x \cdot y) \cdot z & = & x \cdot (y \cdot z)  \\
1 \cdot x & = & x          \\
x  \cdot (y + z) & = & x \cdot y + x \cdot z\\
(x + y)  \cdot z & = & x \cdot z + y \cdot z
\end{eqnarray}

{\bf end}
\newline

These axioms generate a wealth of properties of $+, -, \cdot$ with which we will assume the reader is familiar. We will write $x-y$ as an abbreviation of $x + (-y)$. We notice that $x \cdot 1 = x$ is not implied by these axioms, whereas $x \cdot 0 = 0$ and $0 \cdot x = 0$ are derivable: $0 = x \cdot 0  - x \cdot 0= x \cdot (0 + 0) - x \cdot 0 = (x \cdot 0 + x \cdot 0) - x \cdot 0 = x \cdot 0 + (x \cdot 0 - x \cdot 0) = x \cdot 0$, the other proof is similar.

\subsection{Some concepts of ring theory}\label{AboutRings}
A ring is  {\it commutative} if it satisfies :
\[\forall x \exists y. (x \cdot y  = y \cdot x).\]
A ring is called {\it von Neumann regular} (regular for short) if it satisfies:
\[\forall x \exists y. (x \cdot y \cdot x = x).\]
An element $y$ with $x \cdot y \cdot x = x$ is called a {\it pseudoinverse} of $x$. 
Moreover, an element $y$ with $y \cdot x = 1$ is called an {\it inverse} of $x$. Indeed every inverse is a pseudoinverse as well.

Following e.g., \cite{Ja90}: a ring is {\it strongly regular}  if it satisfies:
\[\forall x \exists y. (x \cdot x \cdot y = x).\]
A regular ring is called {\em unit regular} if it satisfies:
\[\forall x \exists y \exists z. (x \cdot y \cdot x = x\, \&\, y \cdot z = 1).\]

An {\it idempotent} is an element $e$ of a ring that satisfies $e \cdot e = e$. An element is 
$c$ is {\it central} if it
satisfies $\forall x.(x \cdot c = c \cdot x)$.

There is an equivalent definition of strong regularity that is closer to our objectives:  a ring is {\em strongly regular} if it is regular and its idempotents are central.

An element $x$ of a ring is {\it nilpotent} if some power of it equals 0, i.e., 
\[\exists n. x^n = 0.\] 
A ring is {\em reduced} if it has no non-zero nilpotent elements, i.e., 
\[\forall  x. (x \cdot x = 0 \Longrightarrow x=0).\] 
It is immediate that strongly regular rings are reduced. One can prove as a corollary of Proposition \ref{PIL} below, that a reduced regular ring is strongly regular.  

In \cite{AC48} it is shown (using the first definition) that a strongly regular ring is regular. According to
\cite{SAT69} in a strongly regular ring every idempotent is central. Therefore both definitions coincide.  Below we will need this information, but we will provide complete proofs and we will indicate how in an alternative and shorter exposition use could have been made from this equivalence.

We notice that the central elements of a ring constitute a subring. Further, in commutative rings regularity and strong regularity coincide. 

\subsection{Fields and skew fields}
A {\it skew field}, also called a {\it division ring}, is a unital ring that satisifies the general inverse law ($Gil$):
\[\forall x \neq 0. \exists y. (y \cdot x = 1).\]
A division ring is a {\it division algebra} if it is finitely generated over its centralizer subring. The quaternions as designed by Hamilton constitute a division algebra.
A commutative skew field is called a {\it field}.

In skew fields, 1 is a right unit as well as a left unit, and left inverses are also right inverses. Every skew field is a strongly regular ring. Indeed consider $x$, then if $x = 0$, we have 
$x = x \cdot x \cdot y$ (for any $y$) and if $x \neq 0$ then there is  some $y$ with $x \cdot y = 1$ for which of course $x = x \cdot x \cdot y$. Skew fields are also unit regular rings. Consider again $x$: if $x=0$ then $x = x \cdot 1 \cdot 1$ and $1 \cdot 1 = 1$, while with $x \neq 0$ and $y$ its inverse: 
$x \cdot y \cdot x = x \cdot 1 = 1$ and $y \cdot x = 1$.

\subsection{The intended meaning of meadows and skew meadows}
Meadows and skew meadows are concepts created for the following purpose: meadows are supposed to be models of the equational theory of zero totalized fields and skew meadows are supposed to be models of the equational theory of zero totalized skew fields (i.e., skew fields with zero totalized division).  Thus, they have an intended meaning that is independent of their axiomatic definition.

It is due to Theorem \ref{COMP}, that the intended meaning of skew meadows is captured by the technical definition that we will provided below. The corresponding fact was established for meadows in \cite{BHT07}.

The intended meaning of meadows and skew meadows implies that if one considers structures which violate equations about division true in division algebras these algebras will not be called meadows. A typical example of such an equation is $x \cdot x^{-1} = x^{-1} \cdot x$. An algebra that fails to comply with this equation cannot be a meadow or a skew meadow.

\section{The signature and axioms of skew meadows}
To the signature  $\Sigma_{RU}$ of rings with unit, we add an inverse operator  $^{-1}$ to form the primary signature 
$\Sigma_{Md}$, which we will use for totalized division rings and meadows:\\

{\bf signature}  $\Sigma_{Md}$

{\bf import}  $\Sigma_{RU}$

{\bf operations}

$ ^{-1}  \colon ring \to ring$

{\bf end}
\newline 

As we insist that all operations are interpreted as total functions $x^{-1}$ must be defined for all $x$. A division operator which has been made total is called a totalized division operator. Fields and skew fields can be enriched to $\Sigma_{Md}$ algebras by extending the signature and defining $0^{-1}= a$ for some appropriate $a$. A (skew) field thus obtained is called {\em $a$-totalized}. In the sequel we will only consider the case $a=0$ and work with zero totalized fields and skew fields.

A {\it skew meadow} is a $\Sigma_{RU}$ algebra which satisfies $RU$ and in addition these two equations:
\newline

{\bf equations}   {\it SkMd}

{\bf import} {\it RU}, $\Sigma_{Md}$
\begin{eqnarray}
(x^{-1})^{-1}& = & x \label{Ref} \\
x \cdot (x \cdot x^{-1}) & = & x \label{Ril}
\end{eqnarray}

{\bf end}
\newline

Axiom \ref{Ref} is called \axref\, for {\it reflection}. Axiom \ref{Ril} is called \ril\, for {\it restricted inverse law}. Together these axioms imply $0^{-1} = 0$ because 
$0^{-1} = 0^{-1} \cdot 0^{-1} \cdot  (0^{-1})^{-1} = 0^{-1} \cdot 0^{-1}  \cdot 0 = 0$.

A {\em meadow} is a commutative skew meadow.  In terms of axioms, the equations for skew meadows result from the equations for meadows, by simply dropping commutativity of multiplication and including a second distributivity law; this is the same modification that is needed to move from unital commutative rings to arbitrary rings with a left unit.

A zero totalized field is a meadow and a zero totalized skew field is a skew meadow. Just as meadows capture the equational theory of fields with zero totalized division, in a similar fashion, skew meadows have been designed to capture the equational theory of skew fields with zero totalized division. 

A skew meadow is (an enrichment of) a strongly regular ring. A meadow is an enrichment of a commutative regular ring (which must be strongly regular as well). In \cite{BHT07} it was shown that every commutative regular ring can be expanded with an inverse operator to a meadow. Below, in Proposition \ref{SRRtoSkMd}, we will find that in the noncommutative case every strongly regular ring can be expanded to a skew meadow. 

The non-commutative case makes distinctions between axioms that are useful equivalents in the commutative case. For example, the equation
\[\forall x.(x\cdot x^{-1}\cdot x=x)\] 
is called  $\pil$ and cannot be assumed to be equivalent with $\ril$.

\subsection{Multiplicative semigroups}\label{MultSGs}
Forgetting 0, + and - a meadow restricts to an inverse semigroup, which is a special case of regular semigroups. Thus a meadow combines a commutative ring and an inverse semigroup and a skew meadow combines a possibly noncommutative ring with an inverse semigroup. The theory of inverse semigroups provides information about the existence of an inverse operation in regular semigroups. From that perspective the information in our Propositions \ref{Unique} and \ref{SRRtoSkMd} below that strongly regular rings can be uniquely equipped with an inverse operator that satisfies both laws is known as a result in semigroup theory. We have included the results and proofs for completeness of our exposition.

Inverse semigroups have been introduced to formalize properties of function spaces. More specifically, inverse semigroups are to partial symmetries what groups are to symmetries \cite{Law98}. We have not found evidence that inverse semigroups have been used to formalize totalized division in the mathematical literature.

\section{Theory of skew meadows}
In this section we provide a number of logical and structural results concerning skew meadows. We start with some preparations.
\subsection{Auxiliary operators}
When working with meadows we will frequently use some auxiliary notation in order to improve readability or conciseness. We formalize the use of these inessential notations by including them in an extended specification
{\it Aux}. Of course, these operators can always be removed from any specification or proof in favor of their explicit definitions.
\newline

{\bf signature}  $\Sigma_{Aux}$

{\bf sorts}     $ring$

{\bf operations}  

$-^{2} \colon  ring \to ring$;

$\overline{-}   \colon  ring \to ring$;

$ 1_{-} \colon  ring \times ring \to ring$;

$ Z \colon  ring \to ring$;

$ \frac{-}{-} \colon  ring \times ring\to ring$

{\bf end}\\
\newline

{\bf equations}   {\it Aux}

{\bf import} $\Sigma_{Md},\,\Sigma_{Aux}$

\begin{eqnarray}
x^2 &=& x \cdot x\\
\overline{x} &=& x^{-1}\\
1_x &=&  x \cdot \overline{x}\\
Z(x) &=& 1 - 1_x\\
\frac{x}{y} &=& x \cdot \overline{y}
\end{eqnarray}

{\bf end}
\newline

\subsection{Some properties of skew meadows}
In this section we will state and prove a number of useful facts valid in the variety of skew meadows.
\begin{propos} $x \cdot x = 0 \rightarrow x  = 0$.
\end{propos} 
\begin{proof}
If $x \cdot x = 0$ then also $x \cdot x \cdot x^{-1} = 0 \cdot x^{-1} = 0 $. Now using  \ril\, one obtains $x=0$.
\end{proof}

\begin{propos} $x^{-1} \cdot x^{-1} \cdot x = x^{-1}$.
\end{propos} 
\begin{proof}
Combine \ril\, and \axref.
\end{proof}

\begin{propos} $x \cdot y = 1 \rightarrow x = y^{-1}$.
\end{propos} 
\begin{proof}
Assume $x \cdot y = 1$, then $1 = x \cdot y = x \cdot y \cdot y \cdot y^{-1} = 1 \cdot y \cdot y^{-1} =
y \cdot y^{-1} $. Therefore $x = x \cdot 1 =  x \cdot (y \cdot y^{-1}) = (x \cdot y) \cdot y^{-1} =
 1 \cdot y^{-1} =  y^{-1}$. 
\end{proof}
An obvious and useful consequence of this proposition is $1^{-1} = 1$.
\begin{propos} $x\cdot x^{-1}$ is an idempotent: $(x\cdot x^{-1})(x\cdot x^{-1})=x\cdot x^{-1}$.
\end{propos} 
\begin{proof}
We have:  $ x\cdot x^{-1}=(x\cdot x\cdot  x^{-1})\cdot  x^{-1}$ and therefore
$( x\cdot x^{-1})^2=
(x\cdot x\cdot  x^{-1}\cdot x^{-1})^2=
x \cdot x \cdot (x^{-1} \cdot x^{-1} \cdot x)\cdot x \cdot x^{-1}\cdot x^{-1}=
(x \cdot x \cdot x^{-1}) \cdot x \cdot x^{-1} \cdot x^{-1}=(x \cdot x \cdot x^{-1}) \cdot x^{-1}=
x\cdot x^{-1}$.
\end{proof}

\begin{propos}\label{PIL}

(i) $\forall x.(x \cdot 1 = x)$; 

(ii) $\forall x.(x\cdot x^{-1}\cdot x=x)$ (i.e., \pil); and,

(iii)  every idempotent is central: if  $e\cdot e=e$ then  $\forall x(e \cdot x=x\cdot e)$.
\end{propos} 
\begin{proof} This proof uses the following idea taken from  \cite{G79} (chapter 3).
For every $y$ and every idempotent $e$ we have:
$[e\cdot y\cdot (1-e)]^2=0$ as  $(1-e)\cdot e=0$,
and for the same reason also $[(1-e)\cdot y\cdot e]^2=0$.
Therefore we have $e\cdot y\cdot (1-e)=0$ and $(1-e)\cdot y\cdot e=0$. 

Using distributivity we find
 $ e\cdot y\cdot 1-e\cdot y\cdot e=0$  and $ 1 \cdot y\cdot e-e\cdot y\cdot e=0$.
Combining these two equations we have: $ e\cdot y\cdot 1= 1\cdot y\cdot e$ and then 
$ e\cdot y\cdot 1= y\cdot e$.

In particular with $e = y \cdot y^{-1}$ we have $y \cdot y^{-1} \cdot y \cdot 1= y \cdot y \cdot y^{-1}$, and using \ril\, we find $y \cdot y^{-1} \cdot y\cdot 1= y$. Right-multiplying both sides with $y \cdot y^{-1}$ we have $y \cdot y^{-1} \cdot y\cdot 1 \cdot y \cdot y^{-1}= y \cdot y \cdot y^{-1}$ and then 
$y \cdot y^{-1} \cdot y\cdot (1 \cdot y \cdot y^{-1})= y$. Now using $RU$ we have 
$y \cdot y^{-1} \cdot y\cdot (y \cdot y^{-1})= y$ which can be written as 
$y \cdot y^{-1} \cdot (y\cdot y \cdot y^{-1})= y$. Now \ril\, gives 
$y \cdot y^{-1} \cdot y= y$ (i.e., \pil). Combining this fact  with the equation
$y \cdot y^{-1} \cdot y\cdot 1= y\cdot y \cdot y^{-1}$ that was established above one finds 
$y \cdot 1 = y$ for arbitrary $y$. Now using the equation $ e\cdot y\cdot 1= y\cdot e$ which has already been established, we have $ e\cdot y= y\cdot e$ as desired.
\end{proof}

This proof makes use only of the fact that the meadow is reduced. An alternative (and shorter) proof using the literature cited in Section \ref{AboutRings} is as follows: the meadow is a strongly regular ring which must therefore be regular (which proves \pil) and moreover its idempotents are central. Because 1 is an idempotent it is central and therefore $x \cdot 1 = 1 \cdot x = x$. 

\begin{propos} $x \cdot x^{-1} = x^{-1} \cdot x$.
\end{propos} 
\begin{proof}
$x \cdot x^{-1} = x \cdot (x^{-1} \cdot x^{-1} \cdot x) = (x \cdot x^{-1} \cdot x^{-1}) \cdot x =
[(x \cdot x^{-1}) \cdot x^{-1}] \cdot x =$\\
$[x^{-1} \cdot (x \cdot x^{-1})] \cdot x =  x^{-1} \cdot x $.
\end{proof}

\begin{cor} Every skew meadow is a Dedekind finite ring, i.e., it satisfies:\\
$x \cdot y = 1 \rightarrow y \cdot x = 1$.
\end{cor}
\begin{proof} If $x \cdot y = 1$ then $y = x^{-1}$ and thus $y \cdot x = x^{-1} \cdot x = x \cdot x^{-1} = x \cdot y = 1$.
\end{proof}

\begin{propos}\label{expURR} Skew meadows are expansions of unit regular rings.
\end{propos} 
\begin{proof}
Given $x$, if $1_x = 1$ then $x^{-1}$ is a unit. Otherwise, let $y = (1-1_x) + x^{-1}$ and 
$y^{\prime} = (1-1_x) + x$. Now $x \cdot y \cdot x = x \cdot 1 \cdot x - x \cdot 1_x \cdot x +
x \cdot x^{-1} \cdot x =  x \cdot x - x \cdot x + x = x$ and moreover: $ y \cdot y^{\prime} =
[(1-1_x) + x^{-1}] \cdot [(1-1_x) + x] = 
(1-1_x) \cdot (1-1_x) + (1-1_x) \cdot x + x^{-1} \cdot  (1-1_x) + x^{-1} \cdot x =
1- 1_x - 1_x +1_x \cdot 1_x + x - 1_x \cdot x + x^{-1}  - x^{-1} \cdot  1_x + x^{-1} \cdot x =
1- 1_x - 1_x +1_x  + x - x + x^{-1}  - x^{-1} + 1_x =1.$
\end{proof}

\begin{propos}\label{Unique}
The inverse function is uniquely determined in a strong sense: For a given $x$ if
there is 
some $y$ such that:
($x \cdot y  \cdot x=x$   or  $x  \cdot x \cdot y=x$)
and
($y \cdot x \cdot y=y$  or  $y \cdot y \cdot x=y$),
then $y=x^{-1}$.
\end{propos} 
\begin{proof}
If $x \cdot y \cdot x=x$   then $x \cdot y$ is idempotent so that $x \cdot x \cdot y=x$ also holds.
The same holds for $y \cdot x \cdot y=y$. Therefore we have in any case that $x \cdot x \cdot y=x$ and $yyx=y$.

We denote $e= x\cdot x^{-1}$ and make use the fact that it is an idempotent and that it
commutes with every element,  
that it is equal to $x^{-1}x$, and that $x=ex$:

$ y=y \cdot y \cdot x=y \cdot y \cdot (e \cdot e \cdot e \cdot x)=e \cdot e \cdot y \cdot e \cdot y \cdot x= x^{-1} \cdot x^{-1} \cdot x \cdot x \cdot y\cdot e \cdot y \cdot x=x^{-1} \cdot x^{-1} \cdot x\cdot e \cdot y \cdot x
=x^{-1} \cdot x^{-1} \cdot x^{-1}  \cdot x \cdot x \cdot y \cdot x=
x^{-1} \cdot x^{-1} \cdot x^{-1} \cdot x \cdot x=x^{-1} \cdot x^{-1} \cdot x=x^{-1}.$
\end{proof}

\begin{propos} (i) $x \cdot 1_x= x = 1_x \cdot x$\\
(ii) $1_{x \cdot y} \cdot 1_x \cdot 1_y = 1_{x \cdot y}$\\
(iii) $1_x \cdot 1_y = x \cdot y \cdot \overline{y} \cdot \overline{x}$\\
(iv) $\overline{x \cdot y} = \overline{y} \cdot \overline{x}$.
 
\end{propos} 
\begin{proof}
(i) is just \ril, (and for the second equality \pil).\\
(ii)$1_{x \cdot y} \cdot 1_x \cdot 1_y = x \cdot y \cdot \overline{x \cdot y} \cdot 1_x \cdot 1_y =
x \cdot 1_x \cdot y \cdot 1_y  \cdot \overline{x \cdot y} =
x \cdot y  \cdot \overline{x \cdot y} =1_{x \cdot y}.$\\
(iii) immediate.\\
(iv) $\overline{x \cdot y} =
 \overline{x \cdot y} \cdot \overline{x \cdot y} \cdot x \cdot y=
 \overline{x \cdot y}^2 \cdot x \cdot x \cdot \overline{x} \cdot y \cdot 1_y = 
  \overline{x \cdot y}^2 \cdot x  \cdot 1_y \cdot x \cdot \overline{x} \cdot y =
 \overline{x \cdot y}^2 \cdot x \cdot y \cdot \overline{y} \cdot x \cdot \overline{x} \cdot y  =
  \overline{x \cdot y} \cdot \overline{y} \cdot x \cdot \overline{x} \cdot y  =
  \overline{x \cdot y} \cdot \overline{y} \cdot x \cdot 1_x \cdot \overline{x} \cdot y \cdot 1_y  =
  \overline{x \cdot y}  \cdot 1_x \cdot 1_y  \cdot \overline{y} \cdot x \cdot \overline{x} \cdot y  = 
     \overline{x \cdot y} \cdot x \cdot y  \cdot \overline{y} \cdot \overline{x}  \cdot \overline{y} \cdot x \cdot \overline{x} \cdot y  =
     \overline{x \cdot y} \cdot x \cdot y  \cdot \overline{y} \cdot \overline{x}  \cdot \overline{y} \cdot 1_x \cdot y  = 
     \overline{x \cdot y} \cdot x \cdot y  \cdot 1_x \cdot \overline{y} \cdot \overline{x}  \cdot \overline{y}  \cdot y  =
     \overline{x \cdot y} \cdot (x \cdot y)  \cdot 1_x \cdot \overline{y} \cdot \overline{x}  \cdot 1_y  =
 \overline{x \cdot y} \cdot (x \cdot y)  \cdot 1_x \cdot 1_y \cdot \overline{y} \cdot \overline{x}   = 
 \overline{x \cdot y} \cdot (x \cdot y)^2  \cdot \overline{y}\cdot \overline{x} \cdot \overline{y} \cdot \overline{x}   =   
 x \cdot y  \cdot \overline{y}\cdot \overline{x} \cdot \overline{y} \cdot \overline{x}   =
 1_x \cdot 1_y \cdot \overline{y} \cdot \overline{x}   =   1_x \cdot \overline{y} \cdot \overline{x} =    
 \overline{y}\cdot 1_x \cdot \overline{x} = \overline{y} \cdot \overline{x}
$.
\end{proof}
\begin{propos} $1_{x \cdot y} = 1_{y \cdot x}$.
\end{propos} 
\begin{proof}
$ 1_{x \cdot y} = x \cdot y \cdot \overline{x \cdot y} = x \cdot y \cdot \overline{y} \cdot \overline{x} =
x \cdot 1_y \cdot \overline{x} = x \cdot \overline{x} \cdot 1_y = 1_x \cdot 1_y = 1_y \cdot 1_x =
1_{y \cdot x}
$.
\end{proof}

\begin{cor} Every skew meadow  satisfies:\\
$x \cdot y = 0 \rightarrow y \cdot x = 0$.
\end{cor}
\begin{proof} If $x \cdot y = 0$ then 
$y \cdot x = y \cdot x \cdot 1_{y \cdot x} = y \cdot x \cdot 1_{x \cdot y} = 
y \cdot x \cdot x \cdot y \cdot \overline {x \cdot y} = 
y \cdot x \cdot 0 \cdot \overline {x \cdot y} = 
0$.
\end{proof}

\subsection{Embedding theorem and completeness results}
Because the axioms for meadows are equations, and as they hold in all skew fields, \emph{all subalgebras of products of skew fields are skew meadows.} Generalizing the results obtained in \cite{BHT07} to the noncommutative case, we show the converse:

\begin{theorem}\label{EMB}
Every non-trivial skew meadow is isomorphic to a $\Sigma_{Md}$ substructure of a
product of zero totalized skew fields.
\end{theorem}

\begin{proof}
The plan of the proof is as follows. Let $R$ be a  skew meadow in which 0 and 1 are different. For
every element $x\neq 0$ there is a homomorphism (with respect to
all the operations, including the inverse) of the structure $R$
into a zero totalized skew field $R_x$, which takes $x$ away from 0. These
homomorphisms can then be combined into a monomorphism of $R$ into
the product of the skew fields $R_x$ fo all non-zero $x$.

The proof is done in three steps: given an element $x\in R$ we do
the following:

 1. Define a ring homomorphism $h$ from $R$ onto the ring
 $e_x\cdot R$ with its new unit $1_x$, and show that $x$ is mapped
 to an invertible element.

 2. Prove that a skew meadow modulo a maximal ideal is a ring with
 division. If $p$ is the projection of the skew meadow
 $1_x\cdot R$ onto this ring with division then $p\cdot h$ is a ring homomorphism
onto a skew field that maps $x$ away from 0.

 3. We conclude the proof by showing that a ring homomorphism from the reduct of a skew meadow
into a skew field is a meadow homomorphism when the inverse is appropriately defined on the skew field.

\medskip

From now on an {\it ideal} in a ring will mean a two sided (left
and right) ideal.
\medskip

 \emph{Step 1}:
 Because $1_x$ is an idempotent in the meadow $R$ it is central, and for that reason 
$1_x\cdot R$ is a (two sided) ideal, $1_x$ is a unit in this
 ideal, $x$ and $x^{-1}$ are in this ideal, inverse to each other, and
 $h(z)=1_x\cdot z$ is a ring homomorphism. In $1_x\cdot R$ the interpretation of inverse as the
 mapping of $1_x\cdot z$ to $e_x\cdot z^{-1}$ makes the ring $1_x\cdot
 R$ into a meadow. The proof is easy,
 using the axioms of a skew meadow and the properties in the propositions above.

\medskip

 \emph{Step 2}:  Let $S$ be a skew meadow and $J$ a maximal (two sided)
ideal. Then
 $S/J$ is a division ring, and (trivially) every invertible element of $S$
 has a non zero image.

\medskip

$J$ is a maximal ideal. Assume that $x$ is not in $J$ and look at
the subset $J+1_x\cdot S$, which is also a (two sided ideal).
Therefore it is all of $S$, and for some $r\in S$ and $j\in J$ we
have $j+1_x\cdot r=1$. Since $1_x=x\cdot x^{-1}$ (and also $
x^{-1}\cdot x$) we conclude that in the quotient $[x]\cdot
[x^{-1}\cdot r]= [1]$.
 Therefore every non zero element in  $S/J$ has an inverse.

\medskip

 \emph{Step 3}: If $H$ is a ring homomorphism from a skew meadow into a
skew field then $H$ preserves also inverses. Therefore it is a
meadow homomorphism.

\medskip

 If $H(x)=0$ then $H(1_x)=H(x\cdot x^{-1})=H(x)\cdot
H(x^{-1})=0$ so that also $H(x^{-1})=H(1_x\cdot
x^{-1})=H(1_x)\cdot H(x^{-1})=0=(H(x))^{-1}$. We assume therefore
that $H(x)\neq 0$. Then $H(x)=H(1_x\cdot x)=H(1_x) \cdot H(x)$ which
proves that $H(1_x)=1$, by cancellation in skew fields. In other words
$1=H(x\cdot x^{-1})=H(x) \cdot H(x^{-1})$, which proves that
$H(x^{-1})=(H(x))^{-1}$.

\end{proof}

\subsection{Strongly regular rings can be turned into meadows}
By Proposition \ref{PIL}  every skew meadow is an expansion of a strongly regular ring. Our next
observation is the
converse: In every strongly regular ring there is a unique way to choose for each
element one of its
inverse elements and define an inverse function which is also reflexive:

\begin{theorem}\label{SRRtoSkMd}
If $R$ is a strongly regular ring then there is an inverse function (necessarily unique
by Proposition \ref{Unique}) that turns it into a skew meadow.
\end{theorem}
\begin{proof}
We will show first that for all the elements $y$ that are pseudo inverse to $x$, 
i.e., $xyx=x$ the idempotent $xy$ is the same element 
and the element $xyx$ is the same element. We wil then show that
with the definition   $x^{-1}=yxy$ we obtain a skew meadow. We proceed in four steps.

\medskip

a) If $y$ and $y'$ are pseudoinverses of $x$ then  $x\cdot y=x\cdot y'$.
Indeed, first notice that $x \cdot y$ is idempotent: $(x \cdot y) \cdot (x \cdot y) = (x \cdot y \cdot x) \cdot y = x \cdot y$. Then we calculate as follows:

$x\cdot y'= (x \cdot y \cdot x)\cdot y'= (x \cdot y) \cdot (x\cdot y') = (x \cdot y')\cdot (x \cdot y) =
 (x \cdot y' \cdot x) \cdot y=x \cdot y$.

\medskip

If $y$ is a pseudo  inverse of $x$ then we call the product $x\cdot y$
(which is independent of the choice of $y$)  the {\it local unit}
of $x$, and we denote it by $1_x^l$. Trivially $1_x^l\cdot x=x\cdot
1_x^l$ and easily $1_x^l\cdot 1_x^l=1_x^l$.

\medskip

b) For every pseudo  inverse $y$ of $x$, $1_x^l\cdot y$ is also a pseudo inverse
 of $x$. Indeed:
$x \cdot (1_x^l\cdot y) \cdot x=(x \cdot 1_x^l)\cdot y \cdot x=x \cdot y \cdot
x=x$. Moreover, for every pseudo inverse $y$ the product $1_x^l\cdot y$ yields
the same pseudo inverse element. First notice that for a pseudo inverse $y$ of $y \cdot x$ is an idempotent: $(y \cdot x) \cdot (y \cdot x)= y \cdot (x \cdot y \cdot x) =  y \cdot x$. Suppose that $u$, $y$ and $y'$ are pseudo inverses of $x$ then
$1_x^l\cdot y=(x \cdot u) \cdot y= y \cdot (x \cdot u) = y \cdot (x \cdot y' \cdot x) \cdot u=
 (y \cdot x) \cdot (y' \cdot x) \cdot u=
 (y' \cdot x) \cdot (y \cdot x) \cdot u= y' \cdot (x  \cdot y \cdot x) \cdot u=
 y' \cdot (x \cdot u)=
 (x \cdot u) \cdot y'= 1_x\cdot y$.

\medskip

Let $y$ be any pseudo inverse for $x$ we define $x^{-1}$ as follows: $x^{-1}=1_x^l\cdot y$. It has just been demonstrated that this definition is
independent from the choice of $y$. With $y$ a pseudo inverse for $x$ we find $1_x = x \cdot x^{-1} = x \cdot 1_x^l \cdot y = x \cdot (x \cdot y) \cdot y = (x \cdot y) \cdot x \cdot y = x \cdot y = 1_x^l$.

\medskip

c) By (b) $x^{-1}$ is a pseudo inverse of $x$, and for all pseudo inverses $y$
we have $1_x\cdot y=x^{-1}$ therefore:
 $1_x\cdot x^{-1}=x^{-1}$. Moreover we know that $x^{-1}\cdot (x^{-1})^{-1}$ is and idempotent.
 We use this to show that $1_x=1_{x^{-1}}$ as follows: 
 $1_{x^{-1}}=x^{-1}\cdot (x^{-1})^{-1}= 1_x\cdot x^{-1}\cdot (x^{-1})^{-1}=
 x\cdot x^{-1}\cdot (x^{-1}\cdot (x^{-1})^{-1})=
 x\cdot (x^{-1}\cdot (x^{-1})^{-1})\cdot x^{-1}=
 x\cdot (x^{-1}\cdot (x^{-1})^{-1}\cdot x^{-1})= 
 x\cdot x^{-1}=1_x $.
 
\medskip 
 
d) After these preparations we can show that this inverse operation turns the ring into
a skew meadow. \ril\, is trivial since $x^{-1}$ is a pseudo
inverse of $x$, and it remains to show that \axref\, also holds and therefore
$(x^{-1})^{-1}=x$. Indeed:
$ (x^{-1})^{-1}=(x^{-1})^{-1}\cdot 1_{x^{-1}}=(x^{-1})^{-1}\cdot
1_x=(x^{-1})^{-1}\cdot x^{-1}\cdot x=1_{x^{-1}}\cdot x=1_x\cdot
x=x$.
\end{proof}
 
\subsection{A completeness result}\label{COMP}
As in \cite{BHT07}, a completeness result follows from the embedding theorem.
 
\begin{theorem}\label{comp} Every equation valid in all skew fields with zero totalized division is true in the variety of skew meadows and for that reason derivable from   {\it SkMd}.
\end{theorem}
\begin{proof} A trivial skew meadow satisfies all equations and for that reason all equations valid in all skew fields. Now consider a non-trivial skew meadow. According to Theorem \ref{EMB} it is embedded in a subalgebra of a product of skew fields. Now equations true of all skew fields are true in products of skew fields and in all subalgebras of such products including the given non-trivial skew meadow. 
\end{proof}

Another interesting consequence of the embedding theorem is that \emph{a non-trivial skew meadow must be infinite}. According to Wedderburn's Small Theorem, non-commutative skew fields must have characteristic 0. A product of algebras with characteristic 0 has characteristic 0 as well, which implies that its minimal subalgebra is infinite.

\section{Inversion rings}
Inversion rings are like skew meadows but without the restriction that idempotents are central. Because idempotents are central in all skew fields, inversion rings with non-central idempotents will not be referred to as meadows. The ring that underlies an inversion ring must be regular but need not be strongly regular. 

Because strong regularity is implied by the axioms for skew meadows some weakening of the axioms needs to take place.  Instead of \ril\, we will use a (potentially) weaker axiom which we have called \pil\, for {\it pseudo inverse law}. In Proposition \ref{PIL} it was noticed that skew meadows satisfy \pil. The connection between the various notions is now as follows: skew meadows are inversion rings with central idempotents and meadows are commutative skew meadows.

The phrase inversion ring can be easily adapted to fields as follows: an {\em inversion field} is a zero totalized field and a {\em skew inversion field} is a zero totalized skew field. We notice that an inversion skew field is an inversion ring in which all idempotents are central and where the only idempotents are 0 and 1. 

A generalization of an inversion skew field arises as follows: a {\it semi-inversion skew field} is an inversion skew field in which the only central idempotents are 0 and 1. This definition gives rise to the question: {\em Can the embedding theorem and the completeness theorem be generalized to semi-inversion skew fields}. Stated in other words: Is the equational theory of semi-inversion skew fields finitely based?
\newline

\subsection{Axioms for inversion rings}\label{IR}
These axioms are weaker than the axioms for skew meadows.
\newline

{\bf equations}   {\it IR}

{\bf import} {\it RU}, $\Sigma_{Md}$
\begin{eqnarray}
x \cdot 1 &=& x\\
(-x)^{-1} &=& -(x^{-1})\\
(x^{-1})^{-1}& = & x \\
x \cdot (x^{-1} \cdot x) & = & x 
\end{eqnarray}

{\bf end}
\newline

An {\em inversion ring} is a model of  {\it IR}.  
It is immediate that $x \cdot x^{-1}$ is an idempotent in an inversion ring. It is also the case that $0^{-1} = 0$ in any inversion ring. We have included 
$(-x)^{-1} = -(x^{-1})$ because  this equation expresses a very important symmetry. It is easily derivable from the axioms of skew meadows. Currently we have no proof that it is not derivable from the other axioms  of inversion rings.  The same holds for $x \cdot 1 = x$.

\subsection{Pseudo-commutative inversion rings}
An inversion ring is {\em pseudo-commutative} if it satisfies:
$\forall x \forall y.[(x \cdot y)^{-1} = y^{-1} \cdot x^{-1}].$
After removing a redundant axiom one obtains the following axiomatization of pseudo-commutative inversion rings.
\newline

{\bf equations}   {\it PCIR}, $\Sigma_{Md}$

{\bf import} {\it RU}
\begin{eqnarray}
(x \cdot y)^{-1} &=& y^{-1} \cdot x^{-1}\\
(x^{-1})^{-1} &=& x\\
x \cdot (x^{-1} \cdot x) & = & x 
\end{eqnarray}

{\bf end}
\newline

A ring will be called {\em distinctly regular} if its satisfies the following property:
\[\forall x \exists! y. (x \cdot y \cdot x = x\, \& \, y \cdot x \cdot y = y).\]
It is immediate that $x \cdot x = x \rightarrow x = x^{-1}$ in any distinctly regular ring. Strongly regular rings are distinctly regular as a corollary to Proposition  \ref{SRRtoSkMd}. The importance of distinct regularity for inversion rings is implied by the next Proposition.
\begin{propos}\label{DRtoIR}A distinctly regular ring can be expanded to an inversion ring.
\end{propos}
\begin{proof} Obviously on defines $x^{-1} $ as the unique $y$ such that $x \cdot y \cdot x = x$ and 
$ y \cdot x \cdot y = y$. \ril\, is immediately satisfied. We will prove the other two axioms. Consider $u = x^{-1} \cdot 1$. Now $x \cdot u \cdot x = x \cdot x^{-1} \cdot (1 \cdot x) = x \cdot x^{-1} \cdot x = x$ and $ u \cdot x \cdot u = 
x^{-1} \cdot (1 \cdot x) \cdot x^{-1} \cdot 1 = x^{-1} \cdot x \cdot x^{-1} \cdot 1 = x^{-1}  \cdot 1 = u$. Thus $u = x^{-1}$ for all $x$. Now substituting $x^{-1}$ for $x$ we find $x = (x^{-1})^{-1} = u^{-1} = u^{-1} \cdot 1 = (x^{-1})^{-1} \cdot 1 = x \cdot 1$.

Next consider $u = (-x)^{-1}$. We have $ (-x) \cdot u \cdot (-x) = -x$ and $u \cdot (-x) \cdot u = u$ from this we find
$ x \cdot (-u) \cdot x = x$ and $-u \cdot x \cdot -u = -u$. Distinct regularity implies $-u = x^{-1}$ which implies 
$u = -x^{-1}$ which is the required fact.

We notice that it follows from distinct regularity that the expansion is unique.
\end{proof}

\begin{propos}\label{ICtoDR} A regular ring in which all idempotents commute is distinctly regular.
\end{propos}
\begin{proof} The result and its proof are valid in semigroups  but we phrase both in terms of rings.
let $R$ be a regular
ring in which idempotents commute. Then every element of S has at least
one inverse. Suppose that $a$ in $S$ has two pseudoinverses $b$ and $c$, i.e.,
$a\cdot b \cdot a = a, b \cdot a \cdot b = b, a \cdot c \cdot a = a$ and $c \cdot a \cdot c = c$.
Then $a \cdot b,\, a \cdot c.\, b \cdot a$ and $c \cdot a$ are idempotents as immediate consequences of
 these assumptions. 

Now $b = b \cdot a \cdot b = b \cdot (a \cdot c \cdot a) \cdot b = b \cdot a \cdot c \cdot (a \cdot c) \cdot (a \cdot b) = b \cdot a \cdot c \cdot (a \cdot b) \cdot (a \cdot c) = (c \cdot a) \cdot (b \cdot a) \cdot b \cdot a \cdot c = c \cdot a \cdot b \cdot a \cdot c = c \cdot a \cdot c = c.$
 \end{proof}
 
 \begin{propos}\label{DRtoIC} In a distinctly regular  all idempotents commute.
\end{propos}
\begin{proof} 
This result is known in semigroup theory (see \cite{How95} Thm. 5.1.1) but we include a proof to make the paper self contained.

We assume that for every element $x$ in the ring there is a
unique $y$ such that $x\cdot y \cdot x=y$ and $y \cdot x \cdot y=x$. This $y$ is called the
inverse of $x$ and is denoted by $x^{-1}$. We notice that:

\bigskip

1. If $y=x^{-1}$ then $x=y^{-1}$, since the same pair of equations
testify to both.

\bigskip

2. If $e$ is an idempotent then $e^{-1}=e$, since $e=e \cdot e \cdot e$ is the
only condition involved.

\bigskip

{\bf Claim 1:} Let $e$ be an idempotent. If $e \cdot x=x$ then
$x^{-1} \cdot e=x^{-1}$. If $x \cdot e=x$ then $e \cdot x^{-1}=x^{-1}$.

\medskip

{\bf Proof:} $x \cdot (x^{-1} \cdot e) \cdot x=x \cdot x^{-1} \cdot (e \cdot x)=x \cdot x^{-1} \cdot x=x$. On the other
hand $x^{-1} \cdot e\cdot x\cdot x^{-1} \cdot e=[x^{-1} \cdot (e\cdot x)\cdot
x^{-1}] \cdot e=[x^{-1}\cdot x\cdot x^{-1}] \cdot e=x^{-1} \cdot e$. Therefore
$x^{-1} \cdot e$ is inverse to $x$ and by uniqueness it equals $x^{-1}$.
The other direction is similar.

\bigskip

{\bf Claim 2:} If $e$ and $f$ are idempotents then so is also
$e \cdot f$.

\medskip

{\bf Proof:} Let $I=(e \cdot f)^{-1}$.  We have $e\cdot e \cdot f=e \cdot f$ so that
$I \cdot e=I$, by lemma 1. Similarly $f \cdot I=I$. Therefore $I=I\cdot e \cdot f\cdot
I=(I \cdot e) \cdot (f \cdot I)=I \cdot I$. Hence $I$ is idempotent so that $I=I^{-1}=e \cdot f$ and
in particular $e \cdot f$ is idempotent.

\bigskip

Now we can prove that any two idempotents commute.
Let $e$ and $f$ be idempotents. Then $e \cdot f$ is an idempotent and
$(e \cdot f)^{-1}=e \cdot f$. We show that $ \cdot e$ is also inverse to $e \cdot f$ so that
$e \cdot f=f \cdot e$ by the uniqueness of inverse. Indeed
$(e \cdot f) \cdot (f \cdot e) \cdot (e \cdot f)=e \cdot ( \cdot f) \cdot (e \cdot e) \cdot f=
e \cdot f \cdot e \cdot f=e \cdot f$, as $e$ $f$ and $e \cdot f$ are
idempotents. Similarly $(f \cdot e) \cdot (e \cdot f) \cdot (f \cdot e)=f \cdot e$.

 \end{proof}
 
\begin{propos}\label{RRCI2PCIR} A regular ring with commuting idempotents can be expanded to a pseudo-commutative inversion ring.
\end{propos}
\begin{proof} Consider elements $x$ and $y$. We have $x \cdot y \cdot (y^{-1} \cdot x^{-1}) \cdot x \cdot y = 
x \cdot (y \cdot y^{-1}) \cdot (x^{-1} \cdot x) \cdot y = x \cdot (x^{-1} \cdot x) \cdot (y \cdot y^{-1}) \cdot y = x \cdot y$ and similarly $(y^{-1} \cdot x^{-1}) \cdot x \cdot y \cdot (y^{-1} \cdot x^{-1}) = y^{-1} \cdot x^{-1}$. In view of 
Proposition \ref{ICtoDR} the ring is distinctly regular and we can infer that $(x \cdot y)^{-1} = y^{-1} \cdot x^{-1}$.
\end{proof}

\begin{propos}\label{selfinv} In a pseudo-commutative inversion ring where all idempotents $e$ satisfy $e = e^{-1}$ products of idempotents are idempotent.
\end{propos}
\begin{proof} Consider idempotents $e$ and $f$. We have $e \cdot f =  e \cdot f \cdot (e \cdot f)^{-1} \cdot e\cdot f = 
e \cdot f \cdot f^{-1} \cdot e^{-1} \cdot e \cdot f = e \cdot f \cdot f \cdot e \cdot e \cdot f = e \cdot f \cdot e \cdot f $.
\end{proof}

\begin{propos} A distinctly regular inversion ring  is pseudo-commutative.
\end{propos}
\begin{proof}
Using Proposition \ref{DRtoIC} the ring has commuting idempotents. Now using \ref{RRCI2PCIR}
it must be pseudo-commutative as the expansion with an inverse operator must be unique. 

\end{proof}

Several questions remain unanswered: {\em Is there a finite equational specification of the class of directly regular inversion rings?  Is there a finite equational specification of the class of regular inversion rings in which all idempotents are equal to their own inverse?}

\subsection{Inversion compatible rings}

We call a ring {\em inversion compatible} if it can be expanded with an inverse operator into an inversion ring. Proposition \ref{DRtoIR} implies that distinctly regular rings are inversion compatible. Obviously every inversion ring is an expansion of a regular ring. Trivially every inversion compatible ring can be expanded to an inversion ring. We have not been able to answer the following question: {\em Are all regular rings inversion compatible?} We expect this not to be the case.

We will now examine a particular example in meticulous detail. It proves the existence of a non strongly regular and non-pseudo-commutative inversion ring. It can be concluded that {\em the replacement of \ril\, by \pil\, constitutes a weakening of the axioms}.

We have not yet developed any structure theory for non-central inversion rings. All technical work on skew meadows depends on centrality and non-central inversion rings have to be investigated from scratch.

\subsection{Matrix rings over a field of characteristic different from 2}
The $2 \times 2$ matrices over zero-totalized rational numbers $\mathrm{M}_2(\rat_0) $ constitute a non-commutative ring. Although familiar to all mathematicians we will spell out the details of this matter because of the importance of this example. A matrix $X$ has the form 
\[ X = 
\left(\begin{array}{cc}
		x_{11} & x_{12}\\
		x_{21} & x_{22}
       \end{array}
\right).
\]
Let, in addition to $X$, the matrix $Y$ be given by:
\[ Y = 
\left(\begin{array}{cc}
		y_{11} & y_{12}\\
		y_{21} & y_{22}
       \end{array}
\right).
\]
Then we have the two constants for unital rings:
\[ 
0 = 
\left(\begin{array}{cc}
		0 & 0\\
		0 & 0
       \end{array}
\right),\,
1 = 
\left(\begin{array}{cc}
		1 & 0\\
		0 & 1
       \end{array}
\right)
\]
and operators
\[ 
-X = 
\left(\begin{array}{cc}
		-x_{11} & -x_{12}\\
		-x_{21} & -x_{22}
       \end{array}
\right),\,
X + Y = 
\left(\begin{array}{cc}
		x_{11} + y_{11} & x_{12} + y_{12}\\
		x_{21} + y_{21} & x_{21} + y_{22}
       \end{array}
\right)
\]
and
\[
X \cdot Y = 
\left(\begin{array}{cc}
		x_{11} \cdot y_{11} + x_{12} \cdot y_{21} & x_{11} \cdot y_{12} + x_{12} \cdot y_{22} \\
		x_{21} \cdot y_{11} + x_{22} \cdot y_{21} & x_{21} \cdot y_{12} + x_{22} \cdot y_{22} 
		\end{array}
\right) 
\]
An important observation is that
\[ 
\left(\begin{array}{cc}
		0 & 0\\
		1 & 0
       \end{array}
\right) \cdot
\left(\begin{array}{cc}
		0 & 0\\
		1 & 0
       \end{array}
\right) = 
\left(\begin{array}{cc}
		0 & 0\\
		0 & 0
       \end{array}
\right).
\]
Clearly $x \cdot x = 0 \rightarrow x = 0$ fails in matrix rings. Because \ril\, implies $x \cdot x = 0 \rightarrow x = 0$, irrespectively of how an inverse is defined no expansion of the matrix ring to a skew meadow is possible. There might however be an expansion possible to a very skew meadow. It will now be established that this is indeed the case.

\subsection{Expanding the matrix ring with an inverse operator}
The matrix $X$ is called regular if $Det(X) = x_{11} \cdot x_{22} - x_{12} \cdot x_{21} \neq 0$. We abbreviate $D= Det(X)$ and if $X$ is regular it is invertible and an explicit formula for the inverse exists:
\[
X^{-1} = 
\left(\begin{array}{cc}
		\frac{x_{22}}{D} & \frac{-x_{12}}{D}\\
		\frac{-x_{21}}{D} & \frac{x_{11}}{D}
       \end{array}
\right).
\]
In this case we have: $(X^{-1})^{-1} = X$, $X \cdot X^{-1} = X^{-1} \cdot X = 1$ and for that reason also \pil: $X \cdot (X^{-1} \cdot X) = X$. 

We will now define a suitable inverse in all other cases as well. 
Nonregular matrices may have either four zeros (the case $X=0$), three zeros (which splits in four cases), two zeros (either in the same column or in the same row), or no zeros at all and the second column equal to a scalar product of the first column. Together these are seven cases which can be dealt with (almost) independently. In the case of four zeros we take $0^{-1} = 0$ which clearly satisfies both \axref and \pil.
In the case of three zeros we have $x \neq 0$ and the two diagonal cases:
\[ 
\left(\begin{array}{cc}
		x & 0\\
		0 & 0
       \end{array}
\right)^{-1} =
\left(\begin{array}{cc}
		x^{-1} & 0\\
		0 & 0
       \end{array}
\right),\, 
\left(\begin{array}{cc}
		0 & 0\\
		0 & x
       \end{array}
\right)^{-1} =
\left(\begin{array}{cc}
		0 & 0\\
		0 & x^{-1}
       \end{array}
\right)
\]
In these cases \ril\, and \pil\, are immediate and moreover $X \cdot X^{-1} = X^{-1} \cdot X$. Then consider both remaining non-diagonal three-zero cases:
\[ 
\left(\begin{array}{cc}
		0 & 0\\
		x & 0
       \end{array}
\right)^{-1} =
\left(\begin{array}{cc}
		0 & x^{-1}\\
		0 & 0
       \end{array}
\right),\, 
\left(\begin{array}{cc}
		0 & x\\
		0 & 0
       \end{array}
\right)^{-1} =
\left(\begin{array}{cc}
		0 & 0\\
		x^{-1} & 0
       \end{array}
\right)
\]
Now notice:
\[ 
\left(\begin{array}{cc}
		0 & 0\\
		x & 0
       \end{array}
\right) \cdot
\left(\begin{array}{cc}
		0 & x^{-1}\\
		0 & 0
       \end{array}
\right) = 
\left(\begin{array}{cc}
		0 & 0\\
		0 & 1
       \end{array}
\right),\,\mathrm{and}
\left(\begin{array}{cc}
		0 & x^{-1}\\
		0 & 0
       \end{array}
\right) \cdot
\left(\begin{array}{cc}
		0 & 0\\
		x & 0
       \end{array}
\right) =
\left(\begin{array}{cc}
		1 & 0\\
		0 & 0
       \end{array}
\right)
\]
At this stage it is apparent that the we are not dealing with an Abelian skew meadow because 
\[ 
\left(\begin{array}{cc}
		0 & 0\\
		1 & 0
       \end{array}
\right) \cdot
\left(\begin{array}{cc}
		0 & 1\\
		0 & 0
       \end{array}
\right) \neq
\left(\begin{array}{cc}
		0 & 1\\
		0 & 0
       \end{array}
\right) \cdot
\left(\begin{array}{cc}
		0 & 0\\
		1 & 0
       \end{array}
\right) 
\]
In both non-diagonal cases \axref\, is immediate and in the first case \pil\, follows from:
\[ 
\left(\begin{array}{cc}
		0 & 0\\
		x & 0
       \end{array}
\right) \cdot
\left(\begin{array}{cc}
		1 & 0\\
		0 & 0
       \end{array}
\right) = 
\left(\begin{array}{cc}
		0 & 0\\
		x & 0
       \end{array}
\right) =
\left(\begin{array}{cc}
		0 & 0\\
		0 & 1
       \end{array}
\right) \cdot
\left(\begin{array}{cc}
		0 & 0\\
		x & 0
       \end{array}
\right).
\]
In the second non-diagonal case a similar calculation works for \pil:
\[ 
\left(\left(\begin{array}{cc}
		0 & x\\
		0 & 0
       \end{array}
\right) \cdot
\left(\begin{array}{cc}
		0 & 0\\
		x^{-1} & 0
       \end{array}
\right)\right) \cdot 
\left(\begin{array}{cc}
		0 & x\\
		0 & 0
       \end{array}
\right) =
\left(\begin{array}{cc}
		1 & 0\\
		0 & 0
       \end{array}
\right) \cdot
\left(\begin{array}{cc}
		0 & x\\
		0 & 0
       \end{array}
\right)=
\left(\begin{array}{cc}
		0 & x\\
		0 & 0
       \end{array}
\right).
\]
There are four cases with two zeros, writing 2 for 1+1, which are pairwise inverses. We have $x \neq 0$ and $y \neq 0$.
\[ 
\left(\begin{array}{cc}
		x & y\\
		0 & 0
       \end{array}
\right)^{-1} =
\left(\begin{array}{cc}
		(2 \cdot x)^{-1} & 0\\
		(2 \cdot y)^{-1} & 0
       \end{array}
\right),\, 
\left(\begin{array}{cc}
		x & 0\\
		y & 0
       \end{array}
\right)^{-1} =
\left(\begin{array}{cc}
		(2 \cdot x)^{-1} & (2 \cdot y)^{-1}\\
		0 & 0
       \end{array}
\right).
\]
Again \axref\, is immediate. For \pil\, we check:
\[ 
\left(\left(\begin{array}{cc}
		x & y\\
		0 & 0
       \end{array}
\right)\cdot
\left(\begin{array}{cc}
		(2 \cdot x)^{-1} & 0\\
		(2 \cdot y)^{-1} & 0
       \end{array}
\right)\right)\cdot
\left(\begin{array}{cc}
		x & y\\
		0 & 0
       \end{array}
\right) =
\left(\begin{array}{cc}
		1 & 0\\
		0 & 0
       \end{array}
\right) \cdot
\left(\begin{array}{cc}
		x & y\\
		0 & 0
       \end{array}
\right) =
\left(\begin{array}{cc}
		x & y\\
		0 & 0
       \end{array}
\right), 
\]
and
\[ 
\left(\left(\begin{array}{cc}
		x & 0\\
		y & 0
       \end{array}
\right)\cdot
\left(\begin{array}{cc}
		\frac{1}{2 \cdot x} & \frac{1}{2 \cdot y}\\
		0 & 0
       \end{array}
\right)\right)\cdot
\left(\begin{array}{cc}
		x & 0\\
		y & 0
       \end{array}
\right) =
\left(\begin{array}{cc}
		\frac{1}{2} & \frac{x}{2 \cdot y}\\
		\frac{y}{2 \cdot x} & \frac{1}{2}
       \end{array}
\right) \cdot
\left(\begin{array}{cc}
		x & 0\\
		y & 0
       \end{array}
\right) =
\left(\begin{array}{cc}
		x & 0\\
		y & 0
       \end{array}
\right).
\]
The other two cases with two zeros are these:
\[ 
\left(\begin{array}{cc}
		0 & 0\\
		x & y
       \end{array}
\right)^{-1} =
\left(\begin{array}{cc}
		0 & (2 \cdot x)^{-1} \\
		0 & (2 \cdot y)^{-1} 
       \end{array}
\right),\, 
\left(\begin{array}{cc}
		0 & x\\
		0 & y
       \end{array}
\right)^{-1} =
\left(\begin{array}{cc}
		0 & 0\\
		(2 \cdot x)^{-1} & (2 \cdot y)^{-1}
       \end{array}
\right).
\]
Again \axref\, is immediate. For \pil\, we check:
\[ 
\left(\left(\begin{array}{cc}
		0 & 0\\
		x & y
       \end{array}
\right)\cdot
\left(\begin{array}{cc}
		0 & (2 \cdot x)^{-1} \\
		0 & (2 \cdot y)^{-1} 
       \end{array}
\right)\right)\cdot
\left(\begin{array}{cc}
		0 & 0\\
		x & y
       \end{array}
\right) =
\left(\begin{array}{cc}
		0 & 0\\
		0 & 1
       \end{array}
\right) \cdot
\left(\begin{array}{cc}
		0 & y\\
		x & y
       \end{array}
\right) =
\left(\begin{array}{cc}
		0 & 0\\
		x & y
       \end{array}
\right), 
\]
and
\[ 
\left(\left(\begin{array}{cc}
		0 & x\\
		0 & y
       \end{array}
\right)\cdot
\left(\begin{array}{cc}
		0 & 0\\
		\frac{1}{2 \cdot x} & \frac{1}{2 \cdot y}
       \end{array}
\right)\right)\cdot
\left(\begin{array}{cc}
		0 & x\\
		0 & y
       \end{array}
\right) =
\left(\begin{array}{cc}
		\frac{1}{2} & \frac{x}{2 \cdot y}\\
		\frac{y}{2 \cdot x} & \frac{1}{2}
       \end{array}
\right) \cdot
\left(\begin{array}{cc}
		0 & x\\
		0 & y
       \end{array}
\right) =
\left(\begin{array}{cc}
		0 & x\\
		0 & y
       \end{array}
\right).
\]
Finally the case of a nonregular matrix with all elements non-zero works as follows:
\[ 
\left(\begin{array}{cc}
		x & x \cdot y\\
		x \cdot z & x \cdot y \cdot z
       \end{array}
\right)^{-1} =
\left(\begin{array}{cc}
		\frac{1}{4 \cdot x} & \frac{1}{4 \cdot x} \cdot \frac{1}{z}\\
		\frac{1}{4 \cdot x} \cdot \frac{1}{y} & \frac{1}{4 \cdot x} \cdot \frac{1}{y} \cdot \frac{1} { z}
       \end{array}
\right)
\]
\axref\, follows by means of elementary calculation. For \pil\, we check:
\[ 
\left(\left(\begin{array}{cc}
		x & x \cdot y\\
		x \cdot z & x \cdot y \cdot z
       \end{array}
	\right)\cdot
	\left(\begin{array}{cc}
		\frac{1}{4 \cdot x} & \frac{1}{4 \cdot x} \cdot \frac{1}{z}\\
		\frac{1}{4 \cdot x} \cdot \frac{1}{y} & \frac{1}{4 \cdot x} \cdot \frac{1}{y} \cdot \frac{1} { z}
       \end{array}
\right)\right) \cdot
\left(\begin{array}{cc}
		x & x \cdot y\\
		x \cdot z & x \cdot y \cdot z
       \end{array}
	\right) =
\]
\[ 
	\left(\begin{array}{cc}
		\frac{1}{2} & \frac{1}{2 \cdot z}\\
		\frac{z}{2 } &  \frac{1}{2}
       \end{array}
	\right) \cdot
	\left(\begin{array}{cc}
		x & x \cdot y\\
		x \cdot z & x \cdot y \cdot z
       \end{array}
	\right) =
	\left(\begin{array}{cc}
		x & x \cdot y\\
		x \cdot z & x \cdot y \cdot z
       \end{array}
	\right).
\]
The example constitutes an inversion ring which is not pseudo-commutative. Indeed consider the element 
\[ 
P= \left(\begin{array}{cc}
		1 & 0\\
		1 & 0
       \end{array}
       \right).
\]
$P$ is an idempotent but its inverse is not. In a pseudo-commutative inversion ring inverses of idempotents must be idempotent as well: $e^{-1} = (e \cdot e)^{-1} = e^{-1} \cdot e^{-1}$.

We have yet to find an example of a pseudo-commutative inversion ring which is not a skew meadow.
\subsection{Regular and inverse semigroups}
As we have stated already in Section \ref{MultSGs} a skew meadow can be restricted by forgetting addition and subtraction, thus obtaining a semigroup. Regular semigroups are defined as semigroups which satisfy: $\forall x \exists y. (x \cdot y \cdot x = x)$. Clearly every inversion ring is an enrichment of a regular semigroup. 

A more resticted class of semigroups consists of the inverse semigroups. An inverse semigroup is a semigroup that satisfies:
\[\forall x \exists! y. (x \cdot y \cdot x = x\, \& \, y \cdot x \cdot y = y)\]
Here $\exists! y. \phi$ asserts the  existence of a unique $y$ such that $\phi$. Clearly each inverse semigroup is regular. From Proposition \ref{SRRtoSkMd} we find that strongly regular rings are enrichments (by 0, + and -) of inverse semigroups. Thus any skew meadow is a combination of a reduced ring and an inverse semigroup (and conversely).

If a ring is an enrichment of an inverse semigroup it is obviously inverse compatible. We don't know whether or not the converse is true.

Of course inversion compatibility can be viewed as a property of semigroups rather than as a property of rings. All inverse semigroups are inversion compatible, but the converse fails. Now a regular ring is inversion compatible precisely if its multiplicative semigroup is inversion compatible. 

Inversion compatibility for regular semigroups is probably a more easily accessible topic than inversion compatibility for regular rings.

\section{Specification of three zero totalized division rings}
In this section thee algebraic specifications will be proposed each extending the specification of skew meadows.
\subsection{Rational numbers}
The rational number specification from \cite{BT07} can be weakened and commutativity is not essential because $1_x$ is central in any skew meadow, and that restricted form of commutativity satisfies to prove the initiality result. We refer to \cite{BT95,BT07} for the definitions of initial algebra specifications and initial algebra semantics.
\newline

{\bf  equations} {\it $Zero\, totalized\, rationals$}

{\bf import} {\it SkMd,\,Aux}
\begin{eqnarray}
Z(1 + x^2 + y^2 + z^2 + u^2) &=& 0
\end{eqnarray}

{\bf end}
\newline

We will not repeat the argument of \cite{BT07} proving that an initial algebra specification of zero totalized rationals is obtained. Some comments are in order, however. The main thrust of the proof is to demonstrate that each closed term $t$ is provably equal to a term in a set of canonical forms called the transversal. The transversal consists of $0$, and 
$\underline{k} \cdot (\underline{l})^{-1},-\underline{k} \cdot (\underline{l})^{-1}$ for $k$ and $l$ relatively prime positive natural numbers and with $\underline{n}$ denoting the numeral for $n$, i.e., the sum of $n$ 1's. The proof use induction on the structure of $t$. We will only consider the cases 
$t = r \cdot s$ and $t = r + s$ and ignore negations signs to simplify notation. Suppose that $r = \underline{k} \cdot (\underline{l})^{-1}$ and $s = \underline{m} \cdot (\underline{n})^{-1}$. As $n$ and $m$ are sums of four squares plus 1 it can be shown that $ \underline{l} \cdot (\underline{l})^{-1} =
\underline{n} \cdot (\underline{n})^{-1} = 1$. Further it is easy to prove that for all $m$ and $n$,
$\underline{m} \cdot \underline{n} = \underline{m . n}$. 

Now $r + s =  \underline{k} \cdot (\underline{l})^{-1} + \underline{m} \cdot (\underline{n})^{-1} =
\underline{k} \cdot 1 \cdot (\underline{l})^{-1} + \underline{m} \cdot 1 \cdot (\underline{n})^{-1} = 
\underline{k} \cdot \underline{n} \cdot (\underline{n})^{-1} \cdot (\underline{l})^{-1} +
 \underline{m} \cdot \underline{l} \cdot (\underline{l})^{-1} \cdot (\underline{n})^{-1} = 
 \underline{k} \cdot \underline{n} \cdot (\underline{l} \cdot \underline{n})^{-1} +
 \underline{m} \cdot \underline{l} \cdot (\underline{n}\cdot \underline{l})^{-1} = 
  \underline{k . n} \cdot (\underline{l . n})^{-1} +
 \underline{m . l} \cdot (\underline{n . l})^{-1} = 
  \underline{k . n + m . l} \cdot (\underline{n . l})^{-1} =
  \underline{p' . p}  \cdot (\underline{p'' . p})^{-1} = 
   \underline{p'} \cdot \underline{p} \cdot (\underline{p})^{-1} \cdot (\underline{p''})^{-1} = 
    \underline{p'}  \cdot 1 \cdot (\underline{p''})^{-1} = \underline{p'}  \cdot (\underline{p''})^{-1}
  $. Here $k.n + m . l = p' . p$ and $n . l = p'' . p$ with $p$ the GCD of $k.n + m . l $ and $l.n$.

Next we will consider $ t= r \cdot s$. With induction on $n$ one derives that for all $n$: 
$ \underline{n} \cdot x =  x \cdot \underline{n}$ (i.e., all numerals are central in the initial algebra). Now
$r \cdot s =  \underline{k} \cdot (\underline{l})^{-1} \cdot \underline{m} \cdot (\underline{n})^{-1} =
\underline{k} \cdot \underline{m} \cdot (\underline{l})^{-1} \cdot (\underline{n})^{-1} $, and for terms of this form the previous part of the proof has already established a path towards the canonical form.

One may wonder whether a simpler specification is possible for instance by using the equation 
$Z(1 + x^2) = 0$ instead. Now this will not work because the prime field of characteristic 3 satisfies $Z(1 + x^2) = 0$ but fails to be a homomorphic image of the meadow of rational numbers which is immediate by considering the image of $1 = (1+1+1)/(1+1+1)$   which must equal  0 and 1 simultaneously. (At the time of writing we do not know whether or the equation $Z(1 + x^2 + y^2) = 0$ suffices to specify the rational numbers in this context).
\subsection{Complex rational numbers}
As an exercise we specify the zero totalized complex rationals. This specification is an adaptation of the specification 
presented in \cite{BT07}. $c()$ is the complex conjugate. 
It serves as a unary auxiliary function. Just like the rationals the zero totalized complex rationals constitute a commutative meadow while the required amount of commutativity follows from the axioms for skew meadows already. The adaptation of the proof follows the same lines as in the case of the rational numbers. 
\newline

{\bf  equations} {\it $Zero\, totalized\, complex\, rationals$}

{\bf import} {\it SkMd,\,Aux}

{\bf operations}  

$\bfi \colon  \to ring$;

$c \colon ring \to ring$;

\begin{eqnarray}
\bfi^2 &=& -1\\
c(\bfi) &=&-\bfi\\
c(x^{-1}) &= & c(x)^{-1}\\
c(x + y) &= & c(x) + c(y)\\
c(x \cdot y) &= & c(x) \cdot c(y)\\
1_{c(x)} &=&1_x\\
Z(1 + x \cdot c(x)+ y \cdot c(y)) &=& 0
\end{eqnarray}

{\bf end}
\newline

The specification makes use of $c()$ as an auxiliary operator. This specification suggests an obvious question which was first mentioned in \cite{BT06}: can a specification of the zero totalized complex rational numbers be given that makes no use of any auxiliary function.

\subsection{Zero totalized rational quaternions}
The quaternions are a well-known skew field. The rational quaternions constitute its prime 
sub(skew)field. Its expansion to a skew meadow constitutes the zero totalized quaternions. This is a computable algebra. According to \cite{BT95} every computable data type has an initial algebra specification which may make use of auxiliary functions. Here we will make use of a unary auxiliary function $c$. As it turns out many algebras can be specified by means of an initial algebra specification using a only single unary auxiliary function. We are not aware of any theoretical results that indicate why this is the case and 
for what kind of algebras a the use single unary auxiliary function will not be sufficient for giving an initial algebra specification. 

The auxiliary function $c(-)$ is the conjugate for quaternions. It is a division ring (pseudo) homomorphism which sends \bfi, \bfj, \bfk\, to -\bfi, -\bfj\, and -\bfk\, respectively. 

The following set of equations, together with the axioms on inverse and the equations for rings, specifiy the skew field of rational quaternions with zero totalized division as its initial algebra. Interestingly the sum of four squares in the equation that asserts the existence of inverses is now implicit in the multiplication of a quaternion with its conjugate. We omit the correctness proof because it follows the general pattern as given in \cite{BT06,BT07} for the rationals with out any significant complications.
\newline

{\bf  equations} {\it $Zero\, totalized\, rational\, quaternions$}

{\bf import} {\it SkMd,\,Aux}

{\bf operations}  

$\bfi \colon  \to ring$;

$\bfj \colon  \to ring$;

$\bfk \colon  \to ring$;

$c \colon ring \to ring$;
\begin{eqnarray}
\bfi^2 &=& -1\\
\bfj^2 &=& -1\\
\bfk^2 &=& -1\\
\bfi \cdot \bfj \cdot \bfk &=& -1\\
c(\bfi) &=& -\bfi\\
c(\bfj) &=& -\bfj\\
c(x^{-1}) &= & c(x)^{-1}\\
c(x + y) &= & c(x) + c(y)\\
c(x \cdot y) &= & c(y) \cdot c(x)\\
1_{c(x)} &=&1_x\\
Z(1 + x \cdot c(x)) &=& 0
\end{eqnarray}

{\bf end}
\newline

From these equations one easily proves: $c(0) =0$, $c(1) = 1$, $c(k) = -k$ $\bfj^2 =-1$, 
$\bfi \cdot \bfj = \bfk$, $\bfj \cdot \bfk = \bfi$, $\bfk \cdot \bfi = \bfj$, $\bfj \cdot \bfi = -\bfk$,
$\bfk \cdot \bfj = -\bfi$, $\bfi \cdot \bfk = -\bfj$, and
 $c(-x) = -c(x)$ and many other well-known facts about quaternions. Just like in the case of complex numbers the operation $c$ may be viewed as an auxiliary operator in spite of the fact that it is a very familiar one. In both cases (complex rationals and rational quaternions) the question whether or not a specification can be given without an auxiliary operator is open. 
 
 \section{Concluding remarks}
The generalization of our results on meadows in \cite{BHT07,BHT08} to the noncommutative cases is quite satisfactory. Many issues are left open, however, notably the development of a structure theory for non-central inversion rings. Another line of further work is to specify nonassociative algebras with a zero totalized division operator: zero totalized octonions. Nonassociative multiplication is relevant for the subject of division by zero also if one contemplates alternatives containing some form of infinity value that will serve as a proper inverse of zero.  

It is easy to see that the first order theory of fields is 1-1 reducible to the first order theory of meadows which in its turn is 1-1  reducible to the first order theory of skew meadows. Because the first order theory of fields is undecidable so is the first order theory of meadows and so is the first order theory of skew meadows. 

For the equational theory of meadows at least these questions are currently open: (i) can one successfully perform a Knuth-Bendix completion, and:  (ii) is the equational theory of  meadows decidable. The same questions can be posed for skew meadows. For the axioms of pseudo-commutative inversion rings, as well as the axioms for inversion rings, once more the same questions can be posed, such as what rings are inversion compatible, i.e., can be expanded to an inversion ring?

\end{document}